\documentclass[11pt]{paper}

\usepackage{amssymb,amsmath,enumerate,theorem,epsfig,rotating,graphics,changebar,eepic}\usepackage{amsfonts}

\usepackage{a4,latexsym,parskip}

\usepackage{hyperref}

\hypersetup{pdfauthor=MSAS} \hypersetup{pdftitle=DS}

\newcommand{\PP}{\mathbb{P}}

\newcommand{\C} [1][]{\mathbb{C}^{#1}}

\newcommand{\Q} [1] []{\mathbb{Q}_{#1}}

\newcommand{\N} [1][] {\mathbb{N}_{#1}}

\newcommand{\Z}{\mathbb{Z}}

\newcommand{\NN}{\mathcal{N}}

\newcommand{\qed}{\hfill \ensuremath{\Box}}

\theoremstyle{break} \newtheorem{Theorem}{Theorem}

\newtheorem{Proposition}[Theorem]{Proposition}

\newtheorem{Lemma}[Theorem]{Lemma}

\newtheorem{Corollary}[Theorem]{Corollary}

\newtheorem{Remark}[Theorem]{Remark}

\newenvironment{proof}{{\it Proof.}}{\qed}

\begin{document}

\setlength{\unitlength}{1cm}

\title{On the uniqueness of elliptic K3 surfaces with maximal singular fibre}

\author{Matthias Sch\"utt, Andreas Schweizer}

\date{\today}

\maketitle

%


\begin{abstract}
We explicitly determine the elliptic $K3$ surfaces with section and maximal singular fibre. If the characteristic of the ground field is different from $2$, for each of the two possible maximal fibre types, $I_{19}$ and $I^*_{14}$, the surface is unique. In characteristic $2$ the maximal fibre types are $I_{18}$ and $I^*_{13}$, and there exist two (resp. one) one-parameter families of such surfaces.
\end{abstract}



%
%


\maketitle

\section{Introduction}


The possible singular fibres of complex elliptic surfaces
have been classified by Kodaira  in \cite{Kodaira}.
Tate proved in \cite{Tate} that the same classification persists over any perfect field.
A singular fibre always consists of a finite number of rational curves
arranged in compatibility with an extended Dynkin diagram of ADE type.
For instance, a fibre of type $I_n (n>2)$ in Kodaira's notation consists of a cycle 
of $n$ rational curves,
each meeting its two neighbours transversally.
A singular fibre of an elliptic surface is called maximal
if the number of its components attains the maximum possible
within the specified deformation class of elliptic surfaces.

The problem of maximal singular fibres is classically solved for rational elliptic surfaces
(types $I_9, I_4^*, II^*$).
Shioda in \cite{Sh-max} treated the case of elliptic K3 surfaces with section in characteristic $0$ where the maximal fibres have type $I_{19}$ and~$I_{14}^*$. Using the Artin invariant, the first author proved in \cite{S-max} that these are also the maximal fibre types in characteristic $p>0$ if $p$ is odd. Meanwhile in characteristic $2$, types $I_{18}$ and $I_{13}^*$ were shown to be maximal.

Then we considered in \cite{SS} the maximal fibres of general elliptic surfaces with section over $\PP^1$. We proved that in general the maximal fibres are strictly larger in positive characteristic than in characteristic zero. Moreover, we also derived partial uniqueness results.

In this paper, we prove the uniqueness in the K3 case, a problem suggested to one of us by Shioda.

\begin{Theorem}\label{thm}
Let $p$ be an odd prime and $k$ an algebraically closed field of characteristic $p$. 
For each of the maximal fibre types $I_{19}$ and $I_{14}^*$, there exists an elliptic K3 surface with section over $k$, unique up to isomorphism, with a singular fibre of that type.
\end{Theorem}

For both fibre types we have an obvious candidate at hand: 
the mod $p$ reduction of the corresponding elliptic K3 surface over $\Q$ (cf.~\cite{Sh-max}, \cite{ST}). This approach fails only if $p=2$; 
in this case we will determine the families of elliptic K3 surfaces with section that realise the respective maximal fibres $I_{18}$ and $I_{13}^*$ (Propositions \ref{Prop:13-2}, \ref{Prop:18-2}).

Although we are primarily interested in positive characteristic, our method
at the same time re-proves the maximality and uniqueness in characteristic $0$.

We note that the K3 case is special not only in the sense that the results concerning maximal fibres types from characteristic zero hold. Even over $\C$, the uniqueness breaks down as soon as we consider honestly elliptic surfaces, i.e.~with Euler number $e\geq 36$ (cf.~\cite{Sh-DS}).

For characteristic different from $2$ and $3$ the result on maximality 
and uniqueness can be reformulated in the following elementary way.

\begin{Theorem}\label{thmDS}
Let $k$ be an algebraically closed field of characteristic different from 
$2$ and $3$. Fix $M=3$ or $4$. Let $f, g\in k[T]$ with $\deg(f)=2M$, 
$\deg(g)=3M$. Then $f^3 \neq g^2$ implies $\deg(f^3 -g^2)\ge M+1$. Moreover,
up to affine transformation of $T$ and scaling there exists exactly one 
pair $(f,g)$ with $\deg(f^3 -g^2)=M+1$.
\end{Theorem}

\emph{Proof:} 
If there is no polynomial $a$ with $a^4|f$ and
$a^6|g$, then the equation 
\[
Y^2 =X^3 -3fX-2g
\]
defines an elliptic $K3$ surface with
discriminant $\Delta=-108(f^3 -g^2)$. The fibre at $\infty$ has type $I_n$ if 
$M=4$, and $I^*_n$ if $M=3$, with $n=6M-\deg(\Delta)$. Hence the claim follows from Theorem~\ref{thm}.
On the other hand, if there is a polynomial $a$ with $a^2|f$ and
$a^3|g$, then $a^6|(f^3-g^2)$. In particular, $\deg(f^3-g^2)\geq 6>M+1$.
\qed

The degree of $f^3 - g^2$ (for general $M$) is subject to Hall's conjecture, see \cite{Hall}.
In characteristic zero, this degree is minimal if and only if the associated Belyi function $f^3/g^2$ is branched  above only three points of $\PP^1$ (cf.~\cite{St}). 
Through the $j$-map, this problem connects to elliptic surfaces, see \cite{Sh-DS} for 
characteristic $0$ and \cite{SS} for positive characteristic.
The semi-stable elliptic K3 case where the $j$-map has degree $24$ has been solved  for characteristic zero in \cite{BeuMon}.

Indeed, one can try to prove the uniqueness of the $K3$ surfaces in 
characteristic $p\ge 5$ by spelling out the equation 
$\deg(f^3 -g^2)=M+1$ and solving for the coefficients of $f$ and $g$.
In Sections \ref{s:14,p>2} and \ref{s:19} we will use a different 
Weierstrass equation that makes the calculations much less involved
and that works in characteristic $3$ as well. 
In characteristic $2$ 
a somewhat more structural approach applies, based on results by the 
second author in \cite{Schw}. Here we also make use of a classification 
of wild ramification of singular fibres. Although this might be known 
to the experts, we could not find a reference for it, so we include it 
in Section \ref{s:wild} (Proposition~\ref{Prop:wild}).


\textbf{Conventions:} 
Throughout the paper, an elliptic surface is assumed to have a section and a singular fibre. 
The former condition is mild, since we can always consider the jacobian surface.
The latter rules out products.

For a configuration of singular fibres $I_{n_1},\hdots,I_{n_r}^*$, 
we are going to use the shorthand notation $[n_1,\hdots,n_r^*]$.

\section{The configurations in odd characteristic}

We call an elliptic surface $S$ (in)separable if its $j$-map is (in)separable.
Associated to a non-isotrivial elliptic surface $S$, 
there is a unique separable elliptic surface $X\to\PP^1$ 
such that $S$ arises from $X$ by purely inseparable base change of degree $p^d$.
This is the inseparability degree of the $j$-invariant of $S$.
By construction 
the conductors of $S$ and $X$ have the same degree. 
Hence we can apply the following result to $X$ 
in order to reveal additional information about $S$.


\begin{Theorem}[Pesenti-Szpiro {\cite[special case of Theorem~0.1]{P-S}}]
\label{Thm:P-S}
Let $S\to{\mathbb P}^1$ be a non-isotrivial elliptic surface with
conductor $\NN$ and Euler number $e(S)$. Then
\[
e(S)\leq 6 \,p^d (\text{deg}\, \NN -2).
\]
where $p^d$ is the inseparability degree of the $j$-invariant of $S$.
\end{Theorem}

\begin{Remark}
In \cite{S-max} the Artin invariant \cite{Artin} was used to determine
the maximal singular fibres of elliptic K3 surfaces.
We mention that the inequality from Theorem \ref{Thm:P-S}
allows an approach 
without the Artin invariant. We briefly sketch the main idea.
\par

First of all, if $S$ has a fibre of type $I_n$, then necessarily $n\le 21$ since
there are no elliptic surfaces with conductor of degree
less than $4$. 
In the non-isotrivial case, this is a consequence of Theorem \ref{Thm:P-S} since the Euler number of an elliptic surface equals the degree of the discriminant divisor which is a positive multiple of $12$.
Alternatively one can argue with the Shioda-Tate formula for the Picard number \cite[Corollary 5.3]{ShMW}.
Analogously, $m\le 16$ for fibre type $I^*_m$ (at least in
characteristic different from $2$).

\par

If there is a fibre of type $I_{21}$, $I_{20}$, $I^*_{16}$ or $I^*_{15}$ and
the characteristic is different from $2, 3, 5, 7$, the surface is
separable. Moreover, the degree of the conductor is at most $5$,
which by the inequality of Pesenti and Szpiro contradicts $e(S)=24$.
In small characteristic some more fine-tuning is required.
\end{Remark}

The remaining results of this section would follow immediately from the 
explicit determination of the elliptic surfaces in the next two sections. 
However, since the calculations of the equations are long, in particular
for a fibre of type $I_{19}$, we have decided to also include the 
following structural proofs.

\begin{Lemma}\label{Lem:14-3}
Let $S$ be an elliptic $K3$-surface in characteristic $p\ge 0$ with 
a fibre of type $I^*_{14}$. If $p\neq 2,7$ the configuration of $S$ is 
[14*,1,1,1,1].
\end{Lemma}

\emph{Proof:} 
Since $p\neq 2$, the $j$-invariant has a pole of order $14$, so $S$
is non-isotrivial. Moreover, $S$ is separable, since $p\neq 2,7$.
The degree of the conductor of $S$ is at most $6$. Actually, it has to be
exactly $6$, since otherwise by the bound of Pesenti-Szpiro  
the surface would be rational.

So if the fibres outside $I^*_{14}$ are all multiplicative, the 
configuration must be as stated. To finish the proof, it suffices 
to show that there is no additive fibre outside $I^*_{14}$.

Assume on the contrary that there is another additive fibre. Then 
we can apply a quadratic twist that ramifies exactly at that fibre
and at $I^*_{14}$. The twisted surface has conductor of degree
$5$ while being separable, so by Theorem~\ref{Thm:P-S} it is rational. 
However, the $j$-invariant still has a pole of order $14$, contradiction.
\qed

\begin{Lemma}\label{Lem:19-3}
Let $S$ be an elliptic $K3$-surface in characteristic $p\ge 0$ with 
a fibre of type $I_{19}$. If $p\neq 2,19$ the configuration of $S$ is 
[19,1,1,1,1,1].
\end{Lemma}

\emph{Proof:} 
The proof is similar to that of Lemma \ref{Lem:14-3}. The degree of the 
conductor of $S$ can only be $6$ by Theorem~\ref{Thm:P-S}. If there were an 
additive fibre, we could apply a quadratic twist that ramifies exactly 
at this additive fibre and $I_{19}$. 
The conductor of the resulting 
surface would have degree $7$. 
By Theorem~\ref{Thm:P-S}, the twisted surface is 
rational or $K3$. But this is impossible with a fibre of type $I_{19}^*$.
\qed

\begin{Lemma}
An elliptic $K3$ surface in characteristic $19$ with a fibre of type $I_{19}$
is necessarily inseparable and its configuration must be $[I_{19},II,III]$.
Moreover, the surface is unique. If we place the fibres of type $I_{19}$, $II$ and
$III$ at $\infty$, $0$ and $1$, its equation is
$$Y^2 =X^3 - 3\,T^7\,(T-1)X+2\,T\,(T-1)^{11}.$$
\end{Lemma}

\emph{Proof:} 
If the surface is separable, then by the same proof as for Lemma \ref{Lem:19-3}
the configuration must be $[19,1,1,1,1,1]$. The $j$-invariant of this 
surface gives a map of degree $24$ from ${\mathbb P}^1$ to ${\mathbb P}^1$.
By the Hurwitz formula we obtain
$-1=-24+\frac{1}{2}\deg({\mathfrak D})$ where ${\mathfrak D}$ is the 
different. The points above $0$ have ramification indices that are
divisible by $3$, so they contribute at least $8(3-1)=16$ to the degree
of the different. Similarly the points above $1$ contribute at least
$12(2-1)=12$. Finally the $j$-invariant has a pole of order $19$ which,
since the ramification is wild, contributes at least $19$. Summing up
we obtain the contradiction
\[
46 = \deg(\mathfrak{D}) \geq 16+12+19 = 47.
\]

\par

We have seen that the surface is inseparable. In consequence its configuration
can only be $[19,II,III]$. Hence the surface is the Frobenius base change of a
(rational) elliptic surface with configuration $[1 ,II,III^*]$. Up to isomorphism, this surface is unique (in any characteristic $\neq 2,3$) by \cite[Lemma~8.2]{SS}. It can be given in Weierstrass equation
$$Y^2 =X^3 - 3\,T\,(T-1)^3\,X-2\,T\,(T-1)^5.$$
Frobenius base change gives the claimed equation after minimalising.
\qed



\begin{Remark}
By a similar argument one can easily show that in characteristic $7$
an elliptic $K3$ surface with a fibre of type $I^*_{14}$ is inseparable
and must have configuration $[14^*,II,II]$. Placing these fibres 
at $\infty$, $1$ and $-1$ one obtains the equation 
$$Y^2 =X^3 -3\,(T^2-1)^3\, X- 2\,T^7\,(T^2-1),$$
which is exactly the surface from \cite[Theorem~3.2]{Sh-max}.
\end{Remark}

\section{Fibre type $I_{14}^*$ in characteristic $\neq 2$}
\label{s:14,p>2}

In this section, we prove Theorem~\ref{thm} for the maximal fibre type $I_{14}^*$. In characteristic zero, the uniqueness was proven in \cite{St}.

Let $k$ be an algebraically closed field of characteristic different from $2$. Let $S$ be an elliptic K3 surface over $k$ with a section. By assumption, we can work with an extended Weierstrass form
\begin{eqnarray}\label{eq:Weier}
S:\;\; Y^2 = X^3 + A(T)\, X^2 + B(T)\,X + C(T).
\end{eqnarray}
In this setting, the discriminant is given by
\[
\Delta=-27\,C^2+18\,A\,B\,C+A^2\,B^2-4\,A^3\,C-4\,B^3.
\]
Note that we have four normalisations available: M\"obius transformation in $T$ (e.g.~to fix the images of the singular fibres) as well as rescaling by $(X, Y) \mapsto (\alpha^2 X,  \alpha^3 Y)$. 
Note that over non-algebraically closed fields, this variable change yields the quadratic twist over $k(\sqrt{\alpha})$.
In the following, we will use these normalisations in a convenient way.

\begin{Proposition}\label{Prop:14}
There is a unique elliptic $K3$-surface over $k$ with a fibre of type $I^*_{14}$. Its equation can be given as 
$$Y^2 = X^3 +(T^3 +2\,T)\,X^2 -2\,(T^2+1)\,X +T.$$
This surface has discriminant $\Delta=4T^4 +13T^2 +32.$
\end{Proposition}

\emph{Proof:}
Assume that $S$ has a fibre of type $I_{14}^*$.
We locate the image of the special fibre at $\infty$. If the fibre type is $I_n^*$ with $n\geq 8$, then we can assume after a translation $X\mapsto X+\alpha(T)$ that
\[
A(T)=a_3T^3+a_2T^2+a_1T+a_0,\;\; B(T)=b_2T^2+b_1T+b_0,\;\;C(T)=c_1T+c_0.
\]
Here $a_3\neq 0$, since otherwise $S$ would be rational. 
Hence we can scale such that $a_3=1$, e.g.~via $(X,Y,T)\mapsto (a_3^4X, a_3^6Y,a_3T)$. 
By construction, $\Delta$ has degree at most $10$:
\[
\Delta=\sum_{l=0}^{10} d_l T^l.
\]
We ask for all solutions to the system of equations
\[
d_5=\hdots=d_{10}=0.
\]
In the first instance, 
\[
d_{10}=-4 c_1 + b_2^2,\;\;\;\;
d_9=-4 c_0 + 2 b_1 b_2+2b_2^2a_2-12a_2c_1
\]
give 
\[
c_1 = b_2^2/4,\;\;\; c_0 = b_1b_2/2-a_2b_2^2/4.
\]
Then we claim that $b_2\neq 0$. Otherwise, $d_8=0$ would imply $b_1=0$. Thus $b_0=0$, since $d_6=0$, so $B\equiv C\equiv 0$ which gives a contradiction. 

Since $b_2\neq 0$, we can rescale $(X,Y,T)\mapsto (\beta^6 X,\beta^9Y,\beta^2 T)$ for a root of $\beta^8=(-b_2/2)$ to set $b_2=-2$ while preserving the normalisation $a_3=1$. By a translation in $T$, we then achieve $b_1=0$. This uses up our final normalisation. The vanishing of
\[
d_8=-4(a_1+b_0-a_2^2),\;\;\;
d_7=-4(a_0+a_1a_2+2b_0a_2-2a_2^3)
\]
then gives $a_1$ and $a_0$. It follows that
\[
d_5= -2a_2(b_0-2a_2^2),\;\;\;
d_6= (b_0-2a_2^2+2)\,(b_0-2a_2^2-2).
\]
The vanishing of $d_6$ corresponds to two choices of $b_0$ which give rise to isomorphic elliptic curves. 
To see this, apply the scaling $(a_2, T, x, y) \mapsto (ia_2,  iT,  -ix,  \zeta y)$ with $\zeta^2=i=\sqrt{-1}$. Let $a_2=\mu, b_0=2\mu^2-2$. 
We obtain a one dimensional family parametrising elliptic K3 surfaces with a fibre of type at least $I_{13}^*$:
\begin{eqnarray}
\begin{matrix}
Y^2 &  =  & X^3 + ((2-\mu^2)(\mu+T)+\mu T^2+T^3) X^2\\
&&  + 2 (\mu^2-1-T^2) X -\mu+T.
\end{matrix}
\end{eqnarray}
This family has discriminant
\begin{eqnarray}\label{eq:Delta_14}
\phantom{\Delta\Gamma}
{\small
\begin{matrix}
\Delta & =  & 8\mu T^5 +(4+8 \mu^2) T^4+(32 \mu-16 \mu^3) T^3+(13+24 \mu^2-16 \mu^4) T^2\\
&& \;\;\;\;+(86 \mu-32 \mu^3+8 \mu^5) T+32-35 \mu^2-28 \mu^4+8 \mu^6.
\end{matrix}}
\end{eqnarray}
Hence there is a unique specialisation with a fibre of type $I_{14}^*$ at $\mu=0$. This gives the claimed equation and discriminant for $S$ in any characteristic $p\neq 2$. \qed

\begin{Remark}
It is immediate from the shape of the coefficients $A, B, C$ that $S$ arises from a rational elliptic surface by quadratic base change and twisting. It follows that this rational elliptic surface is uniquely determined by its configuration of singular fibres. The configuration is [7,1,1,III] if $p\neq 7$, and [7,II,III] if $p=7$.
\end{Remark}

\begin{Remark}
The model in Proposition~\ref{Prop:14} relates to the Weierstrass equation in \cite{Sh-max} as follows: The first author exhibited in \cite[\S 6]{S-max}  a quadratic twist of the latter model which has good reduction at $3$. Then use the coordinate $T=\frac 1s$ and twist over $\Q(\sqrt{-1})$.
\end{Remark}

We can also compute the Picard number $\rho(S)$ of the K3 surface from Proposition \ref{Prop:14}.
It turns out that only two cases occur: $\rho(S)=20$ and $22$.
After Artin \cite{Artin}, K3 surfaces with the later property are called supersingular since $\rho(S)=b_2(S)$.

\begin{Corollary}\label{Cor:ss}
Let $S$ be an elliptic K3 surface over $k$ with a fibre of type $I_{14}^*$. Then $S$ is supersingular if and only if $p\equiv 3\mod 4$.
\end{Corollary}

\emph{Proof:} By Proposition~\ref{Prop:14}, $S$ arises from the corresponding elliptic K3 surface $S_0$ over $\Q$ by reduction. Since $S_0$ is a singular K3 surface (i.e.~$\rho(S_0)=20$), it is modular by a result of Livn\'e \cite[Remark~1.6]{L}. The associated newform of weight $3$ has CM by $\Q(\sqrt{-1})$. (In fact, it has level $16$ (cf.~\cite[Table~1]{S-CM}).) 
This implies that the characteristic polynomial of Frobenius on $H_{\acute{e}t}^2(S, {\Q}_\ell)$ has all zeroes of the shape $\zeta p$ for some root of unity $\zeta$ if and only if $p\equiv 3\mod 4$. Hence the Corollary follows from the Tate conjecture \cite{Tate-C} which is known for elliptic K3 surfaces by \cite{ASD}. \qed

\section{Fibre type $I_{19}$ in characteristic $\neq 2$}

\label{s:19}

In this section, we sketch the proof of Theorem~\ref{thm} for the maximal fibre type $I_{19}$ in characteristic $p\neq 2$. In characteristic zero, this was again included in \cite{St}.

We use the same approach and notation as in Sect.~\ref{s:14,p>2}. In particular, the proof first determines the elliptic K3 surfaces with a fibre of type $I_{18}$.
We find two one-dimensional families (cases 1 and 2 in the proof of Proposition \ref{Prop:19}).
The family with the simpler equation is distinguished by the existence of a 3-torsion section.
Then we determine the unique specialisation within the torsion-free family with a fibre of type $I_{19}$.
Since the proof heavily relies on the help of a machine to factor polynomials, we omit some of the details. They can be obtained from the authors upon request.

\begin{Proposition}\label{Prop:19}
The elliptic $K3$-surface over $k$ with a fibre of type $I_{19}$
is unique. Its equation can be given as 
$$
Y^2 = X^3 + (T^4+T^3+3 T^2+1)\, X^2 + 2 (T^3 + T^2 + 2 T)\, X + T^2 + T + 1.
$$
This surface has discriminant $\Delta=4 T^5+5T^4+18T^3+3T^2+14T-31$. 
\end{Proposition}

\emph{Proof:}
Let $S$ be an elliptic K3 surface with a fibre of type $I_{19}$. Let $S$ be given in extended Weierstrass form (\ref{eq:Weier}).
Locating the special fibre at $\infty$, we can assume that
\begin{eqnarray*}
A(T) & = & a_4T^4+a_3T^3+a_2T^2+a_1T+a_0,\\
B(T) & = & b_3 T^3 + b_2T^2+b_1T+b_0,\\
C(T) & = & c_2 T^2 +c_1T+c_0.
\end{eqnarray*}

Since the special fibre is multiplicative, $a_4\neq 0$. 
Since $k$ is algebraically closed, we can scale such that $a_4=1$. In this setting, $\Delta$ has degree at most $14$:
\[
\Delta=\sum_{l=0}^{14} d_l T^l.
\]
We ask for the solutions to the system of equations
\[
d_6=\hdots=d_{14}=0.
\]
In the first instance, we shall ignore $d_6$, thus investigating the special fibre type $I_{18}$.
The vanishing of the polynomials $d_{14}, d_{13}, d_{12}$ determines $c_2, c_1, c_0$. Then 
\[
d_{11}=-b_3^2 a_1+\hdots,\;\;\;d_{10}=-b_3^2 a_0 + \hdots,
\]
so we have to distinguish whether $b_3=0$.

\textbf{1st case}: $\boldsymbol{b_3=0}$. In this case, $d_{11}= b_2(2b_1-b_2a_3)$. The choice $b_2=0$ successively implies $b_1=b_0=0$. In consequence, $B\equiv C\equiv 0$, and equation (\ref{eq:Weier}) becomes singular. Hence we can assume $b_2=1$ after rescaling and obtain $b_1$ from $d_{11}=0$. The successive factorisations of $d_{10},\hdots$ directly give $b_0, a_1$ and $a_0$. Then we normalise by a translation in $T$ to assume $a_3=0$. 
We obtain a family of elliptic K3 surfaces with a fibre of type $I_{18}$ and a 3-torsion point.
\[
Y^2 = X^3 + \left(T^2 + \frac{a_2}{2}\right)^2 X^2 + \left(T^2 + \frac{a_2}{2}\right) X + \frac 14.
\]
The translation $Y\mapsto Y+(T^2+\frac{a_2}2)X+\frac 12$ moves the 3-torsion point to $(0,0)$ so that the family reads
\begin{eqnarray}\label{eq:18,3-torsion}
Y^2 + (2\,T^2 +a_2)\, X\,Y + Y= X^3.
\end{eqnarray}
There are several notable properties of this family:
\begin{itemize}

\item It arises from the (unique) rational elliptic surface with a fibre of type $I_9$,
\begin{eqnarray}\label{eq:9}
E:\;\;\; Y^2 + T\,  X\, Y + Y = X^3
\end{eqnarray}
by the family of quadratic base changes
\begin{eqnarray*}
T\mapsto 2\,T^2+a_2.
\end{eqnarray*}

\item There is no specialisation with a fibre of type $I_{19}$, since on an elliptic K3 surface, the existence of a 3-torsion section predicts that $3\mid n$ for all singular fibres of type $I_n$ with $n>6$.
(Confer \cite{MP} where the quotient by translation by the section is considered.
This quotient is again K3, so the Euler number yields the divisibilty property.)
Another way to deduce the property is to consider the discriminant 
\[
\Delta=(2 T^2 + a_2 -3) (4 T^4 + (6 + 4 a_2) T^2 + 9 + 3 a_2 + a_2^2).
\]

\item It admits a model over $\Q$ with good reduction at $2$: Instead of the base change
$T\mapsto 2T^2 +a_2$ simply apply a base change which is equivalent up to 
M\"obius transformation, e.g.~$T\mapsto T^2+\lambda T$.
\end{itemize}

\textbf{2nd case:} $\boldsymbol{b_3\neq 0.}$ In this case, the vanishing of $d_{11}$ and $d_{10}$ determines $a_1$ and $a_0$. In consequence,
\begin{eqnarray}\label{eq:h}
d_9=-\frac 2{b_3} \underbrace{(b_1b_3-b_2^2+b_2b_3a_3-b_3^2a_2)}_h b_0 + \hdots .
\end{eqnarray}
The coefficient $h$ of $b_0$ does not vanish because otherwise $d_9=\frac 12 b_3^3\neq 0$. 
Hence we obtain $b_0$ from $d_9=0$. 
This leaves us with polynomials $d_6, d_7, d_8$ in five variables where we have still two normalisations left. In the next instance, we note that
\[
d_7-2a_3d_8=\frac 12 b_3\,(3 b_1 b_3+3 b_2^2-6 b_2 b_3 a_3 -2 b_3^2 a_2+3 b_3^2 a_3^2) 
\]
This factorisation provides us with $a_2$. At this point, we want to analyse the vanishing of $d_8$ independently from $d_6$, i.e.~we first make sure that there is a fibre of type $I_{17}$ and then promote it successively to type $I_{18}$ and $I_{19}$. The numerator of $d_8$ is a complicated polynomial. We shall sketch two ways to solve it. The first normalisation is ad-hoc, while the second will use some extra knowledge.

The first normalisation is a linear transformation in $T$ such that $a_3=1, b_1=0$. (This is possible unless char$(k)=3$ and $b_2=0$. In this special case, we find five single solutions to $d_8$ which are not roots of $d_6$.) Then the numerator of $d_8$ is a polynomial of degree $12$ in $b_2, b_3$. Since every summand has degree at least $10$, this polynomial has degree $2$ in the homogenising variable. Hence it can be solved explicitly. We obtain a one-dimensional rational parametrisation of elliptic K3 surfaces with $I_{18}$ fibres. 
No K3 surface in this family has a torsion section other than the zero section.

The second solution is much more efficient. It was motivated by a private correspondence with N.~Elkies who kindly informed us about an explicit 1-parameter family of elliptic K3 surfaces with an $I_{18}$ fibre and generically trivial group of sections which he had found independently. We therefore decided to choose the normalisation in such a way that it would meet his example after a change of variables.

\textbf{Claim:} For any solution, there are two linear transformations $T\mapsto \alpha T+\beta$ such that after rescaling:
\[
a_4=1, \;\;\; b_2=a_3b_3,\;\;\; b_1=2 b_3.
\]

\emph{Proof:} 
After the transformation in $T$, we rescale $(X, Y) \mapsto (\alpha^4X, \alpha^6Y)$ to retain  our first normalisation $a_4=1$. 
In consequence, the new coefficients read
\[
a_3'=\dfrac{a_3+4\beta}\alpha,\;\;b_3'=\dfrac{b_3}{\alpha^5},\;\;b_2'=\dfrac{b_2+3 b_3\beta}{\alpha^6},\;\; b_1'=\dfrac{b_1+2b_2\beta+3b_3\beta^2}{\alpha^7}.
\]
The first requirement
\[
b_2'=a_3'b_3' \Leftrightarrow b_2+3 b_3\beta =  (a_3+4\beta) b_3
\]
gives $\beta=\frac{b_2}{b_3}-a_3$. Then the second condition $b_1'=2b_3'$ implies
\[
\alpha^2 = -\dfrac {h(a_3, b_1, b_2, b_3)}{b_3^2}
\]
with the polynomial $h$ as in (\ref{eq:h}). By assumption, $h\neq 0$, so the claim follows. \qed

Applying one of the above linear transformations to our elliptic surface, we obtain
\[
d_8=b_3^9 (2 a_3+b_3+4) (2 a_3+b_3-4)/16.
\]
It is easily checked that both solutions are identified under the scaling $(T, a_3)\mapsto (-T, -a_3)$ (which interchanges the two choices of $\alpha$). Setting $b_3=2\lambda, a_3=2-\lambda$, we obtain the second family of elliptic K3 surfaces with a fibre of type $I_{18}$, but without a 3-torsion section:
\begin{eqnarray}\label{eq:Delta-18-1}
\begin{matrix}
Y^2  =  X^3 + (T^4 - (\lambda-2) T^3 + 3 T^2 - (2 \lambda-2) T + 1)\, X^2\\ 
 \;\;\;\;\;\;\;\;\;\;\;\;\;\;\;\;\; + 2 \lambda   (T^3 - (\lambda-2) T^2 + 2 T - \lambda+1)\, X + \lambda^2   ((T+1)^2 - \lambda T).
\end{matrix}
\end{eqnarray}


This has discriminant
\begin{eqnarray*}
\Delta & = & \lambda^3[
(4 \lambda-4)\, T^6
-(8 \lambda^2-24 \lambda+12)\, T^5
+(4 \lambda^3-20 \lambda^2+45 \lambda-24)\, T^4\\
&&\;\;\;\;\;
-(30 \lambda^2-76 \lambda+28)\, T^3
+(13 \lambda^3-52 \lambda^2+66 \lambda-24)\, T^2\\
&&\;\;\;\;\;\;
-(46 \lambda^2-72 \lambda+12)\, T
+32 \lambda^3-96 \lambda^2+37 \lambda-4].
\end{eqnarray*}
Hence the elliptic K3 surface with a fibre of type $I_{19}$ is uniquely obtained as the specialisation at $\lambda=1$. This gives the claimed equation and discriminant. \qed

\begin{Remark}
To relate the models in Proposition~\ref{Prop:19} and in \cite{Sh-max}, we go through the twist in \cite{ST}. Here we have to correct the following typo: One summand of $B$ in \cite[\S 2]{ST} ought to be $15t^4$. Then the equation in Proposition~\ref{Prop:19} is obtained from the translation $x\mapsto x+t+t^2$ after the variable change $t\mapsto -T-1$.
\end{Remark}

\begin{Corollary}
Let $S$ be an elliptic K3 surface over $k$ with a fibre of type $I_{19}$. Then $S$ is supersingular if and only if $p\neq 2$ is not a quadratic residue mod $19$ or $p=19$.
\end{Corollary}

\begin{proof} 
As before, $S$ arises from a singular K3 surface $S_0$ over $\Q$ by reduction by Theorem~\ref{thm}. Hence the claim follows as in the proof of Corollary~\ref{Cor:ss}. (See e.g.~\cite[Corollary~2.2 \& \S 3]{ST}.) 
\end{proof}

\begin{Remark}
In case of characteristic zero, the uniqueness of the underlying complex K3 surface $X$ for the elliptic fibrations in Propositions \ref{Prop:14} and \ref{Prop:19} follows from the Torelli Theorem \cite{PS} since $X$ is determined by its transcendental lattice which corresponds to the unique binary quadratic form of discriminant $-19$ (cf.~\cite{SI}). However, it is a non-trivial problem to show that this guarantees the uniqueness of the respective elliptic fibrations. 
\end{Remark}

\section{Wild ramification of singular fibres}

\label{s:wild}

In the absence of wild ramification, 
the type of a (singular) fibre is determined by the vanishing orders of
the $j$-invariant and the discriminant
(cf.~\cite[p.~365, Table 4.1]{Silverman}). 
In characteristics $2$ and $3$, however, wild ramification
makes it necessary to go through Tate's algorithm
to determine the fibre type.
Although there are papers investigating the case of wild ramification in wide generality (e.g.~\cite{MT}), it seems that there is no reference which gives explicit lower bounds for the index of wild ramification.
For convenience and future reference, we decided to include such a list.

\begin{Proposition}\label{Prop:wild}
Let $E$ be an elliptic curve over a complete valuation ring $R$ with residue characteristic $p=2$ or $3$. Let $\pi$ denote a  uniformizer of $R$. Let $w$ denote the index of wild ramification the special fibre of $E$ at $\pi$. Depending on the reduction type, the following table lists whether we always have $w=0$ or in one case $w=1$, or it gives a sharp lower bound for $w$:
$$
\begin{array}{|c||c|c|}
\hline
\text{fibre type} & p=2 & p=3\\
\hline
\hline
I_n (n\geq 0) & 0 & 0\\
\hline
II & \geq 2 & \geq 1\\
\hline
III & \geq 1 & 0\\
\hline
IV & 0 & \geq 1\\
\hline
I_n^* \;
\begin{matrix}
(n=1) \\ (n\neq 1)
\end{matrix}
& 
\begin{matrix}
1 \\
\geq 2
\end{matrix} & 
\begin{matrix} 
0\\ 0
\end{matrix}\\
\hline
IV^* & 0 & \geq 1\\
\hline
III^* & \geq 1 & 0\\
\hline
II^* & \geq 1 & \geq 1 \\
\hline
\end{array}
$$
\end{Proposition}

The fibre type $I_1^*$ in characteristic $2$ is exceptional due to the fact that it admits wild ramification, but only of index $w=1$.
For all other wild fibres, there is no general upper bound on the index of
wild ramification.
In characteristic $2$, this is visible from the elliptic surface with Mordell-Weil rank $2^r$ that Elkies constructed in \cite{E-MW}.
More recently, Gekeler investigated these questions in \cite{Gekeler}. 
This behaviour is in contrast with the situation 
over a given $p$-adic field or number field. 
(Compare \cite[Theorem 10.4]{Silverman}.)

\begin{Corollary}\label{Cor:wild}
In the notation of Proposition~\ref{Prop:wild}, let the special fibre at $\pi$ be additive. Then the vanishing order of the discriminant $\Delta$ at $\pi$ satisfies
\[
v_\pi(\Delta)\geq \begin{cases} 3, & \text{if } p=3,\\ 4, & \text{if } p=2.\end{cases}
\]
\end{Corollary}

Proposition~\ref{Prop:wild} is easily verified with Tate's algorithm \cite{Tate}. 
Here we shall only give the proof in characteristic $2$ (the case considered in the next two sections). 
We work with the general Weierstrass form
\[
y^2 + a_1 x y + a_3 y = x^3 + a_2 x^2 + a_4 x^2 + a_6\;\;\; (a_i\in R).
\]
In characteristic $2$, the discriminant reads
\begin{eqnarray}\label{eq:d-2}
\Delta=a_1^4 (a_1^2a_6+a_1a_3a_4+a_2a_3^2+a_4^2) +a_3^4 + a_1^3 a_3^3.
\end{eqnarray}
Note that the sharpness of the given bounds will follow immediately from the proof.
We will follow the exposition in \cite{Silverman} and  employ the notation 
\[
a_{i,j}=a_i/\pi^j.
\]
We let $k$ denote the algebraic closure of the residue field of $R$.
Note that the argument works over the residue field itself as well, but then one has to take special care of fibre components which are only defined over some extension.

\begin{proof}[Proof of Proposition~\ref{Prop:wild} in characteristic $p=2$:]
Recall that at an additive fibre the vanishing order of $\Delta$ satisfies
\[
v_\pi(\Delta)=1 + \# (\text{connected components of the fibre}) + w.
\]
Starting with reduced fibres (types $II, III, IV$), we have to prove that $\pi^4|\Delta$ (with equality for type $IV$). 
By a change of variables, we can move the singular point to $(0,0)$. 
Hence $\pi \mid a_3, a_4, a_6$. 
The singular fibre is additive if and only if $\pi \mid a_1$. 
By (\ref{eq:d-2}), this implies $\pi^4 \mid \Delta$. 
In fact, we obtain
\begin{eqnarray}\label{eq:dd1}
\Delta = \pi^4\,a_{3,1}^4 + O(\pi^5).
\end{eqnarray}
Fibre type $III$ moreover requires $\pi^2\mid a_6$, so still $v_\pi(\Delta)\geq 4$. 
For fibre type $IV$, there is another condition $\pi\mid (a_2\,a_{3,1}^2+a_{4,1}^2)$.
The two exceptional components of the fibre are encoded in the distinct zeroes of the polynomial
\[
U^2 + a_{3,1}\,U + a_{6,2} \in k[U].
\]
In particular, this implies $v_\pi(a_{3,1})=0$, so $v_\pi(\Delta)=4$ by (\ref{eq:dd1}).
This completes the analysis of reduced fibres.

Assume now that the singular fibre is not reduced.
Following \cite{Silverman}, we can assume that $\pi|a_1, a_2,\; \pi^2|a_3, a_4,\;\pi^4|a_6$.
In particular, we see that $\pi^8|\Delta$.
We distinguish three cases depending on the number of distinct roots of the polynomial
\[
P(U) = U\,(U^2+a_{2,1}\,U+a_{4,2}) \in k[U].
\]

\textbf{Case 1:}
If $P(U)$ has three distinct roots in $k$, then the singular fibre has type $I_0^*$, and $v_\pi(\Delta)\geq 8$, so $w\geq 2$.

\textbf{Case 2:}
If $P(U)$ has a double root and a simple root in $k$, then the fibre has type $I_n^*$ for some $n>0$.
Translate the root to $0$ so that $\pi\mid a_{4,2}$.
Then
\begin{eqnarray}\label{eq:dd2}
\Delta = \pi^8\,a_{3,2}^4 + O(\pi^9).
\end{eqnarray}
The singular fibre has type $I_1^*$ if and only if the polynomial
\[
Q(U) = U^2 + a_{3,2}\,U + a_{6,4}
\]
has distinct roots in $k$.
Equivalently $\pi\nmid a_{3,2}$, so $v_\pi(\Delta)=8$ by (\ref{eq:dd2}) and $w=1$ as claimed.

Generally, the integer $n>0$ in the fibre type $I_n^*$ is determined by the conditions that after further coordinate changes 
\begin{eqnarray}
\label{eq:tate-odd}
\pi^{l+2}|a_4,\;\, \pi^{2l+2}|a_6, \;\, \pi^{l+1}\,\mid\mid a_3 & \text{if $n=2l-1$ is odd},\\
\label{eq:tate-even}
\pi^{l+2}|a_3, \;\, \pi^{2l+3}|a_6,\;\, \pi^{l+2}\,\mid\mid a_4 & \text{if $n=2l$ is even.} \phantom{-1}\;
\end{eqnarray}
Here $\mid\mid$ denotes exact divisibility. 
To prove the claim $w\geq 2$ for $n>1$, we shall use two-step induction: 
To start the induction, we need that for $n=2$, the a priori lowest order term in (\ref{eq:d-2}) is $a_1^4a_4^2$, so $\pi^{10}|\Delta$. 
If $n=3$, then we have to consider the term $a_1^4a_2a_3^2$ in (\ref{eq:d-2}), so $\pi^{11}|\Delta$. 
To complete the induction, we note that $n\to n+2$ increases the $\pi$-divisibility of every summand in (\ref{eq:d-2}) by at least two by inspection of (\ref{eq:tate-odd}), (\ref{eq:tate-even}). 
Hence $w\geq 2$ if $n\geq 2$. 

\textbf{Case 3:}
If $P(U)$ has a triple root, then we translate it to $0$ again so that $\pi\mid a_{2,1}, a_{4,2}$.
For the discriminant, we still have (\ref{eq:dd2}).

The singular fibre has type $IV^*$ if and only if the polynomial $Q(U)$ has distinct roots in $k$.
As before, this yields $v_\pi(\Delta)=8$, so here $w=0$.

Otherwise, we translate the double root to $0$ so that $\pi\mid a_{3,2}, a_{6,4}$. 
Then 
\[
\Delta = \pi^{10}\,a_{1,1}^4\,a_{4,3}^2 + O(\pi^{11}).
\]
The singular fibre has type $III^*$ if and only if $\pi\nmid a_{4,3}$. Equivalently $v_\pi(\Delta)\geq 10$ and thus $w\geq 1$.
Otherwise we have fibre type $II^*$ with $v_\pi(\Delta)\geq 11$ and again $w\geq 1$.
\end{proof}

The proof in characteristic $3$ is similar, but simpler since we can work with an extended Weierstrass form (\ref{eq:Weier}) and $I_n^*$ fibres do not admit wild ramification.
It is left as an exercise to the reader.

An elliptic curve over a field $K$ of characteristic $2$ with $j(E)\neq 0$ 
can be given in normal form
\begin{eqnarray}\label{eq:normal-2}
Y^2 + XY = X^3 +a_2 X^2 +\frac 1{j(E)}.
\end{eqnarray}

Note, however, that this form in general is not integral or minimal.
A twist replaces the coefficient $a_2$ by $a_2 +D$ while preserving $j(E)$.
Such a twist is trivial (i.e. the two curves are isomorphic over $K$)
if $D$ is of the form $\beta^2 +\beta$ with $\beta\in K$.
If $a_2 =0$, the above normal form has multiplicative fibres at all poles
of $j(E)$.

In characteristic $2$, a fibre of type $I^*_{\nu}$ does not imply that the 
$j$-invariant has a pole at the corresponding place. Actually, every twist
with sufficiently wild ramification will produce a fibre of type $I^*_{\nu}$.
We describe this locally:

\begin{Lemma}\label{Lem:twist}
Let $k$ be an algebraically closed field of characteristic $2$. 
Consider the elliptic curve
$$E:\ \ Y^2+XY=X^3+DX^2+\frac{\pi^r Q}{\pi^{6e}}$$
over $k((\pi))$ where $e\in\Z$, $0\le r\le 5$ and 
$Q=\sum\limits_{i=0}^\infty\omega_i\pi^i$ is a unit in $k[[\pi]]$. 
Let 
$$D=\frac{\delta}{\pi^{2d-1}}+\ldots\in\frac{1}{\pi^{2d-1}}k[[\pi]]$$
with $\delta\in k^\times$ and $d>0$.
\par

If $d>e$ then $E$ has a special fibre of type $I^*_{\nu}$ with 
$$\nu=8d-4-6e+r=8d-4-v(j(E))$$ 
and index of wild ramification $\omega=4d-2$.
\end{Lemma}

\emph{Proof:}
Essentially this is proved in \cite[Lemma 2.1]{Schw}, although the type
of the fibre is not stated explicitly there.
\qed

\section{Type $I_{13}^*$ in characteristic $2$}

\label{s:13-2}

The first approach towards maximal singular fibres on elliptic K3 surfaces in characteristic $2$ might be to mimic the construction for odd characteristic:
reduce the maximal surfaces over $\Q$ from Propositions \ref{Prop:14} and \ref{Prop:19} mod $2$.
However, in characteristic $2$, these equations do only define a quasi-elliptic fibration;
i.e.~every fibre is a curve of genus one, but the general fibre is singular (type $II$).
Such surfaces can only exist in characteristic $2$ and $3$ (cf.~\cite{CD} for instance).

The present quasi-elliptic fibrations have only one reducible singular fibre;
namely it has type $I_{16}^*$.
Both fibrations can be transformed to the Weierstrass form
\[
S:\;\; Y^2 = X^3 + T^3\,X^2 + T.
\]
This surface (and some other models of it) has been studied extensively in \cite{DK}.
In other words, the N\'eron-Severi lattice of $S$ accomodates the fibre through the abstract decomposition $\mbox{NS}(S)=U+D_{20}$, but unlike the usual situation in characteristics different from $2, 3$, the corresponding  fibration is only quasi-elliptic.

Returning to fibre type $I_{14}^*$, it was shown in \cite{S-max} that there is no elliptic K3 surface in characteristic $2$ with such a singular fibre. It was immediate from the proof that the surfaces with a fibre of type $I_{13}^*$ come in a family. In the notation of \cite{S-max}, this family was determined by the choice $c=\sqrt{e}$. To correct a typographical error in \cite{S-max}, we note that this is the exact case when $v_0(\Delta)=21$.

In this section we shall re-prove this result. The motivation to do so is twofold. First, we use an argument which is structural and does not require assistance by a machine. Secondly, we will directly derive the explicit equation of the family.

We now recall some results from \cite{Schw} that are helpful for our classification purposes in characteristic $2$. The setup is an algebraically closed field $k$ of characteristic $2$ and an elliptic curve $E$ over $k(T)$ which is Frobenius-minimal (i.e.~$j(E)$ is separable). 

Note that the places of supersingular reduction of $E$ are exactly the zeroes 
of $j(E)$ that are places of good reduction.

\begin{Lemma}[{\cite[Lemma~2.4 (a)]{Schw}}]\label{Lem:Schw-2}
Let char$(k)=2$. Assume that $E$ is Frobenius-minimal with $j(E)\not\in k$. Let $\beta_i$ denote the finite places of bad reduction of $E$ and 
\begin{eqnarray*}
r_i\in\{0,1,2,3\} \;\;\text{ such that} & r_i\equiv -v_{\beta_i}(j(E))\mod 4,\\
s_i\in\{0,1\} \;\;\text{ such that} & s_i\equiv -v_{\beta_i}(j(E))\mod 2.
\end{eqnarray*}
Define
\[
G(T)=\prod (T-\beta_i)^{r_i},\;\;\; H(T)=\prod (T-\beta_i)^{s_i}.
\]
Let $\vartheta$ be a finite place of supersingular reduction of $E$ and expand\linebreak
 $G(T)=\sum \gamma_i(T-\vartheta)^i$. Let $c_G$ denote be the smallest index $i$ with $4\nmid i$ and $\gamma_i\neq 0$. Then $v_\vartheta(j(E))=12e$ for some $e\in\N$ with 
\[
(T-\vartheta)^{3e-1}|H'(T)\;\;\;\text{ and }\;\;\; 3e\leq c_G.
\]
\end{Lemma}

Now we are ready to prove our first main result in characteristic $2$.

\begin{Proposition}\label{Prop:13-2}
In characteristic $2$, any elliptic $K3$-surface with a fibre of type $I^*_{13}$
occurs in the following family:
$$Y^2 + T X Y + Y = X^3 + \lambda T^3 X^2 + \lambda T, \;\; \lambda\in k^*.$$
All these surfaces have $\Delta=T^3 +1$ and $j=\frac{T^{12}}{T^3 +1}$.
\end{Proposition}

\begin{proof}
We locate the special fibre at $\infty$. Proposition~\ref{Prop:wild} prescribes the vanishing order $v_\infty(\Delta)\geq 21$. Hence it follows from Corollary~\ref{Cor:wild} that all other singular fibres are multiplicative. 

In contrast to other characteristics,
the existence of a fibre of type $I_n^*$ ($n>0$) does not rule out the possibility
of $j(E)=0$ in characteristic $2$. 
In that case we move the maximal fibre to $0$ and have a model
$$Y^2 +T^6 Y=X^3 +a_2 X^2 +a_4 X+a_6.$$
As already explained in \cite{S-max}, the Tate algorithm then shows that
$m\le 9$ for fibres $I^*_m$ (under the condition that we have a K3 surface).
Hence $j(E)\neq 0$, and since the fibre at $\infty$ is additive, there 
is at most one finite supersingular place.

If there is no finite supersingular place (i.e.~$j(E)=\frac 1{\Delta(T)}$), then Tate's algorithm \cite{Tate} for the integral model of the twisted normal form (\ref{eq:normal-2})
\[
Y^2 + XY = X^3 + D(T) X^2 + \Delta(T)
\]
shows that the fibre at $\infty$ can maximally have type $I_{12}^*$. (This was the case $a_1=t^2$ in \cite{S-max}, where the maximal fibre had been located at $0$.)

Otherwise, we locate the finite supersingular place at $T=0$ by a M\"obius transformation. It follows from Lemma~\ref{Lem:Schw-2} that either $v_\infty(\Delta)=21$ with three distinct affine roots or the surface is not Frobenius-minimal. In the latter case, we have $j(E)=\frac{T^{12}}{(aT+b)^2}$ with $b\neq 0$. 
But then the associated separable surface $S(E)$ would have $j(S(E))=\frac{T^6}{aT+b}$
and hence there would be an additive singular fibre at $0$. 
(In characteristic $2$ the $j$-invariant  is $a_1^{12}/\Delta$, so the
multiplicity of a good place in the $j$-invariant must be divisible by $12$.)
Since a purely inseparable base change cannot replace an additive fibre by a smooth or multiplicative fibre, 
we obtain a contradiction.

We now consider the first case, $v_\infty(\Delta)=21$. By Lemma~\ref{Lem:Schw-2}, we obtain 
\[
j(E)=\dfrac{T^{12}}{\varepsilon (T^3+c)}
\]
with $\varepsilon, c\in k^*$. Rescaling, we achieve $c=1$. 
The ``untwisted'' form
$$Y^2 +XY=X^3 +\frac{\varepsilon(T^3 +1)}{T^{12}}$$
has fibres $I_9$ at $\infty$, $I_1$ at the third roots of unity, and
also a singular fibre at $0$. Since in characteristic $2$ every twist 
can be build from twists that ramify at only one place, we can change
the fibres at $0$ and at $\infty$ individually. First we apply a twist 
that makes the fibre at $0$ smooth. In \cite[Proposition 5.1]{Schw} it
was shown that this is only possible with $\varepsilon=1$ and the 
following twist
$$Y^2 +XY=X^3 +\frac{1}{T^3}X^2 +\frac{T^3 +1}{T^{12}}.$$
In other words, our K3 surface with fibre of type $I^*_{13}$ must be a twist 
of the extremal elliptic surface $(\ref{eq:9})$, which has configuration [9,1,1,1].

Now we apply a twist that ramifies only at $\infty$. Lemma~\ref{Lem:twist}
tells us that we will get a fibre of type $I^*_{13}$ at $ \infty$ if and only if
this twist is
$$Y^2 +XY=X^3 +\left(\lambda T+\frac{1}{T^3}\right)X^2 +\frac{T^3 +1}{T^{12}}.$$
The corresponding integral model can be minimalised after the translation $(X, Y) \mapsto (X+T, Y+1)$. This gives the claimed equation.
\end{proof}

A similar approach, following \cite[Lemma~2.4 (b)]{Schw}, can be carried
out in characteristic $3$, for fibre type $I^*_{14}$ as well as for $I_{19}$.
We have omitted this, since the proofs in Sections \ref{s:14,p>2} and 
\ref{s:19} settle all cases of odd characteristic with one calculation.

\section{Type $I_{18}$ in characteristic $2$}

\label{s:18-2}

In \cite{S-max}, it was shown that a fibre of type $I_{19}$ is impossible for an elliptic K3 surface $S$ in characteristic $2$. We briefly sketch how this can be proven purely in terms of Lemma~\ref{Lem:Schw-2}:

Assume that $S$ has a fibre of type $I_{19}$ at $\infty$. Hence $S$ is separable, and $\Delta$ has degree five. We first consider the case with an additive fibre, which we place at $0$. Then, by Corollary~\ref{Cor:wild} and Lemma~\ref{Lem:Schw-2}, there is no other additive fibre and no finite supersingular place. Hence $a_1=T^2$, and the argumentation in the proof of Lemma~\ref{Lem:18} rules this case out.

Otherwise, all fibres are multiplicative. Since $\Delta'$ has degree four, it follows from Lemma~\ref{Lem:Schw-2} that there are two distinct supersingular places (and the configuration is [19,1,1,1,1,1], which also follows directly from Theorem~\ref{Thm:P-S}). Normalise so that $\vartheta_1=0, \vartheta_2=1$. Then $\Delta(T)=\varepsilon H(T)$, where by Lemma~\ref{Lem:Schw-2}
\[
H(T)=T^5+aT^4+T^3+cT^2+e,\;\;\;(e\neq 0, a+c+e\neq 0).
\]
We apply the criterion with the index $c_G$ from Lemma~\ref{Lem:Schw-2}: At $\vartheta_1$, it gives $c=0$. But then the expansion of $H(T)$ at $\vartheta_2$ is 
\[
H(T)=(T+1)^5+(a+1)(T+1)^4+(T+1)^3+(T+1)^2+a+e, 
\]
so the criterion gives a contradiction. Hence there is no elliptic K3 surface in characteristic $2$ with a fibre of type $I_{19}$.

We shall now study the next case of fibre type $I_{18}$. In Section \ref{s:19}, we have exhibited two families of elliptic surfaces over $\Q$ with such a fibre. For the surface in (\ref{eq:18,3-torsion}), we already saw that it has good reduction at $2$. In this section, we will prove that also the second family reduces nicely mod $2$ (Remark~\ref{Rem:transf}), and that any elliptic K3 surface in characteristic $2$ with a fibre of type $I_{18}$ is a member of one of these families (cf.~Proposition~\ref{Prop:18-2}). We first determine the possible configurations:

\begin{Lemma}\label{Lem:18}
Let $S$ be an elliptic K3 surface in characteristic $2$ with a fibre of type
$I_{18}$. Then $S$ is necessarily semistable. 
Moreover, if $S$ is inseparable, it must be the Frobenius base change of the rational elliptic surface $E$ in (\ref{eq:9}), i.e. we have 
$$S:\ Y^2 + T^2XY + Y = X^3$$
with configuration {\rm [18,2,2,2]}. On the other hand, if $S$ is 
separable, the configuration is necessarily {\rm [18,1,1,1,1,1,1]}.
\end{Lemma}

\begin{proof}
We locate the fibre of type $I_{18}$ at $\infty$. Then $\Delta$ must have degree $6$.
Corollary~\ref{Cor:wild} shows that there is at most one additive fibre.

\par

We first assume that there is an additive fibre and locate it at $0$. Hence $T|a_1$.
Since $T^4 |\Delta$ there are at most two more singular fibres. We claim 
that $0$ is the only zero of $j(S)$. This can be seen as follows. If $S$ is separable, then the
polynomial $H(T)$ has degree at most $2$. If $S$ is inseparable, then it 
is the Frobenius base change of a separable elliptic surface with one fibre of type $I_9$ at 
$\infty$, one additive fibre at $0$, and possibly one $I_1$ fibre. Again the 
polynomial $H(T)$ for this surface has degree at most $2$. In either 
case there cannot be a zero of $j(E)$ at a good place by Lemma~\ref{Lem:Schw-2}. This is equivalent to the claim. So in the minimal 
model of $S$
$$Y^2+a_1 XY+a_3 Y=X^3+a_2 X^2+a_4 X+a_6$$ 
we can assume $a_1 =T^2$. From the discriminant
$$\Delta=a_1^6a_6+a_1^5a_3a_4+a_1^4a_2a_3^2+a_1^4a_4^2+a_1^3a_3^3+a_3^4$$
we see that $T$ must divide $a_3$. Hence $\Delta$ is congruent  to $a_3^4$ modulo
$T^8$. In particular $\Delta$ cannot have degree $6$. From this contradiction
we conclude the semistability of $S$.

\par

Now, if $S$ is inseparable, its configuration has to be [18,2,2,2] by Theorem~\ref{Thm:P-S}. Hence it is the Frobenius base change of the unique rational elliptic surface $E$ with configuration [9,1,1,1]  in (\ref{eq:9}).

\par

If $S$ is separable, the conductor must have degree at least $6$ by
Theorem~\ref{Thm:P-S}. Hence the configuration can only be 
[18,1,1,1,1,1,1] or [18,2,1,1,1,1]. But for the second configuration, $\deg H(T)=4$, so
Lemma~\ref{Lem:Schw-2} would imply the contradiction $\deg(a_1)\le 1$.
\end{proof}


\begin{Proposition}
\label{Prop:18-2}
Separable elliptic K3 surface in characteristic $2$ with a fibre of type 
$I_{18}$ come in two families:
\begin{eqnarray}
\label{eq:18-2-3'}
Y^2 + (T^2 +T+1) XY + rY =  X^3,\ \ \ (r\neq 0)
\end{eqnarray}
with $\Delta=r^3(T^6 +T^5 +T^3 +T+1+r)$. This family with a 3-torsion section is obtained from the 
extremal elliptic surface (\ref{eq:9}) by a family of quadratic base changes.
\begin{eqnarray}
\label{eq:18-2-1'}
Y^2 + (T^2+T+1) XY + r (T+1) Y   =  X^3+rT^3 X^2+r^3 T,\ \ \ (r\neq 0)
\end{eqnarray}
with $\Delta=r^3(T^6 +T^5 +rT^4 +T^3 +T+1+r)$.
\end{Proposition}

\begin{proof}
Locating the special fibre at $\infty$ and applying Lemma~\ref{Lem:Schw-2}
we obtain
\[
\Delta(T)=\varepsilon H(T),\;\;\; H(T)=T^6+aT^5+bT^4+cT^3+dT^2+eT+f.
\]
We claim that $a\neq 0$. To prove this, assume on the contrary $a=0$, so $\deg H'\leq 2$. It follows from Lemma~\ref{Lem:Schw-2} that there is at most one finite supersingular place, and that this place has multiplicity $1$ as a root of $a_1$. Hence the elliptic surface is either rational ($\deg a_1=1$), or the fibre at $\infty$ has additive type. This gives a contradiction.

Since $a\neq 0$, $H'$ has degree four. By Lemma~\ref{Lem:Schw-2} and the above argument, there are two distinct supersingular places. For convenience we normalise them to be primitive third roots of unity $\varrho, \varrho^2$: 
\[
A(T)=T^2+T+1.
\]
By Lemma~\ref{Lem:Schw-2}, this implies $e=c=a$:
\[
H(T)=T^6+aT^5+bT^4+aT^3+dT^2+aT+f.
\]
Then the index criterion from Lemma~\ref{Lem:Schw-2} at the supersingular places implies 
\[
(1+ a) \varrho +d = (1+a)\varrho +a+d+1=0.
\]
Hence $a=1, d=0$, and
\[
H(T)=T^6+T^5+bT^4+T^3+T+f.
\]
Note that the variable change $T\mapsto T+1$ preserves $A(T)$ and the shape of $H(T)$ after replacing $f$ by $b+f$. The surface
$$Y^2 + XY = X^3 + \frac{\varepsilon H(T)}{A(T)^{12}}$$
has multiplicative fibres at all poles of the $j$-invariant
$\frac{A(T)^{12}}{\varepsilon H(T)}$ as desired, but it also has
singular fibres at the zeroes of $A(T)$. To make these fibres smooth we
apply a twist that only ramifies at $\varrho$ and $\varrho^2$. By 
Lemma~\ref{Lem:twist} this twist can only be of the form
$$Y^2 + XY = X^3 
+\left(\frac{\alpha_3}{(T-\varrho)^3} +\frac{\alpha_1}{T-\varrho}
+\frac{\beta_3}{(T-\varrho^2)^3} +\frac{\beta_1}{T-\varrho^2}\right)X^2
+\frac{\varepsilon H(T)}{A(T)^{12}}.$$
We now consider the following integral model of $S$:
\begin{eqnarray}\label{eq:int'}
Y^2 + A(T)^2 XY = X^3 + A(T) D(T) X^2 + \varepsilon H(T).
\end{eqnarray}
Here, after a variable change $Y\mapsto Y+A(T) \alpha(T) X$, it suffices 
to allow the following twisting polynomials:
\[
D(T)=d_3T^3+d_2T^2+d_1T+d_0.
\]
In order for the fibres at the supersingular places to be smooth, 
we require that the integral model (\ref{eq:int'}) is not minimal at the supersingular places. 
We pursue these issues simultaneously by going through Tate's algorithm \cite{Tate}. 
We apply several translations and, if necessary, substitutions to increase the vanishing orders of all coefficients of the Weierstrass form successively.
For the first variable change, we introduce new parameters $\lambda, \mu$ with
\[
\lambda^2 = \varepsilon (b+f+1),\;\;\;\;\mu^2 = \varepsilon b.
\]
Then we transform $Y\mapsto Y+\mu T + \lambda$ to obtain
$$
\begin{array}{l}
Y^2 + A(T)^2  XY  =  \\
\;\;\; X^3 + A(T) D(T) X^2 +  A(T)^2 (\mu T + \lambda) X + A(T)^2(\varepsilon A(T)+\mu^2).
\end{array}
$$
Next, $Y\mapsto Y+A(T)\mu$ gives
$$
\begin{array}{l}
Y^2 + A(T)^2 XY = \\
\;\;\;
X^3 + A(T) D(T) X^2 +  A(T)^2 (\mu T +\lambda + \mu A(T)) X + \varepsilon A(T)^3.
\end{array}
$$
Here we have to rule out fibre type $I_n^*$ at the supersingular places. In other words, the following two polynomials must have a triple zero each:
\begin{eqnarray*}
P(Z) & = & Z^3+(d_0+\varrho d_1+\varrho^2 d_2+d_3)Z^2+(\varrho\mu+\lambda) Z+\varepsilon,\\
Q(Z) & = & Z^3+(d_0+\varrho^2 d_1+\varrho d_2+d_3)Z^2+(\varrho^2\mu+\lambda) Z+\varepsilon.
\end{eqnarray*}
We distinguish whether these polynomials are equal or not:

$\boldsymbol{P(Z)=Q(Z)}:$ Denoting the triple root by $r$, we obtain
\[
d_1=d_2,\;\; d_0+d_2+d_3=r,\;\;\mu=0,\;\;\lambda=r^2,\;\;\varepsilon=r^3.
\]
Hence the variable change $X\mapsto r A(T)$ gives
$$
\begin{array}{l}
Y^2 + A(T)^2 XY + r A(T)^3 Y = \\
\;\;\;
X^3 + A(T)^2 (T d_3 + d_2+d_3) X^2 + r^2 A(T)^4 (T d_3 + d_2+d_3).
\end{array}
$$
For the next transformation, we let $u^2=d_2, v^2=d_3$. 
Then the translation $Y\mapsto Y+r A(T)^2(Tv+u)$ gives
$$\begin{array}{r}
Y^2 + A(T)^2 XY + r A(T)^3 Y 
=X^3 +A(T)^2 (T d_3 +d_2+d_3) X^2\\
+rA(T)^4 (Tv+u)X+r^2 A(T)^5 (T v + u+v^2).\\
\end{array}$$
The final step for non-minimality requires that $A(T)^6$ divides the coefficient of the constant term. Since this has only degree $11$, it has to be zero:
\[
u=v^2,\;\;\;v=0.
\]
After minimalising, we obtain the family (\ref{eq:18-2-3'}).

$\boldsymbol{P(Z)\neq Q(Z):}$ In this case, the triple roots differ by a third root of unity. Hence, given a choice of $\varrho$ there are unique $r, \omega$ such that $\omega=\varrho$ or $\omega=\varrho^2$ and
\[
P(Z)=(Z+\omega r)^3,\;\;\; Q(Z)=(Z+\omega^2 r)^3.
\]

Since the choices of $\varrho, \varrho^2$ are permuted by the variable change $T\mapsto T+1$, we can assume that $\omega=\varrho^2$. This implies (by adding the coefficients of $P(Z), Q(Z)$) that
\[
d_1+d_2=r,\;\; d_0+d_1+d_3=0,\;\;\mu=r^2,\;\;\lambda=0,\;\;\varepsilon=r^3.
\]
In consequence, we translate $X\mapsto X+rA(T)(T+1)$. This results in
$$\begin{array}{l}
Y^2 + A(T)^2 XY + r A(T)^3(T+1) Y \\
=X^3 +A(T)^2 (T d_3 +d_3+d_2) X^2 +r^2 A(T)^4 (A(T)(Td_3+d_2)+ T d_1 + d_3).
\end{array}$$
We let $d_1=u^2, d_3=v^2$ and translate $Y\mapsto Y + r A(T)^2 (uT+u+v)$. This gives
\begin{eqnarray*}
 Y^2 + A(T)^2 XY + r A(T)^3(T+1) Y 
= X^3 + A(T)^2 (T d_3 + d_3+d_2) X^2\\
\;\;\; + r A(T)^4 (uT+u+v) X + r^2 A(T)^5 (u A(T) + (u+v+v^2)T + v+r).
\end{eqnarray*}
Again non-minimality requires that the final summand is divisible by $A(T)^6$. Hence
\[
v=r,\;\;\; u=v+v^2.
\]
After minimalising, we obtain the family of elliptic K3 surfaces in the parameter $r$ as
$$
\begin{array}{l}
Y^2 + (T^2+T+1) XY + r (T+1) Y =\\
\;\;\;
X^3+r(r^3+ T r +1) X^2+r^2(rT+T+r)X+r^3 (r+1)
\end{array}
$$
with $H(T)= T^6+T^5+rT^4+T^3+T+r+1$. 
Finally the translation $Y\mapsto Y+r(T+1+r)\,X+r^2$ gives (\ref{eq:18-2-1'}). 
\end{proof}

\begin{Remark}\label{Rem:transf}
For the family (\ref{eq:18-2-3'}), we have seen that it can be obtained from the first family in section \ref{s:19} by reduction mod 2. For the second family, start with (\ref{eq:Delta-18-1}) and apply the variable change
\[
(X, Y, \lambda) \mapsto \left(4X, 8Y, r/4\right).
\]
Then the translation $Y\mapsto Y+\frac 12(T^2+T+1)X+\frac 12 r (T+1)$ gives
\[
Y^2 + (T^2+T+1) XY + r (T+1) Y = X^3 - r T (T^2+2) X^2 - 2 r^2 (T^2+1) X - r^3 T.
\]
This reduces mod $2$ to the claimed equation (\ref{eq:18-2-1'}).
\end{Remark}

%
%
%
%
%
%
%
%

%
%
%
%
 \begin{center}
 \begin{tabular}{ll}
 Matthias Sch\"utt & Andreas Schweizer\\
 Institut f\"ur Algebraische Geometrie\;\;\;\;\;\; & Institute of Mathematics, Academia Sinica\\
 Leibniz Universit\"at Hannover & 6F, Astronomy-Mathematics Building\\
 Welfengarten 1 & No. 1, Sec. 4, Roosevelt Road \\
 30167 Hannover & Taipei 10617\\
 Germany & Taiwan\\
 schuett@math.uni-hannover.de & schweizer@math.sinica.edu.tw
 \end{tabular}
 \end{center}


\begin{thebibliography}{99999}







\bibitem{Artin} Artin, M.: \emph{Supersingular $K3$ surfaces,} Ann.~Sci.~\'Ecole~Norm.~Sup.~(4) {\bf 7} (1974), pp.~543--568.



\bibitem{ASD} Artin, M., Swinnerton-Dyer, P.: \emph{The Shafarevich-Tate conjecture for pencils of elliptic curves on $K3$ surfaces}, Invent.~Math.~{\bf 20} (1973), pp.~249--266.



\bibitem{BeuMon} Beukers, F., Montanus, H.:
\emph{Explicit calculation of elliptic fibrations of K3-surfaces and their Belyi-maps,}  
Number theory and polynomials,  33--51, 
London Math.~Soc.~Lecture Note Ser.~{\bf 352}, Cambridge Univ. Press, Cambridge, 2008.





\bibitem{CD}
Cossec, F.~R., Dolgachev, I.~V.:
\emph{Enriques surfaces}. I. 
Progress in Math.~{\bf 76}. Birkh\"auser (1989).

\bibitem{DK}
Dolgachev, I., Kond\=o, S.:
\emph{A supersingular K3 surface in characteristic 2 and the Leech lattice},
Int.~Math.~Res.~Not.~2003, no.~1, pp.~1--23. 

\bibitem{E-MW}
Elkies, N.~D.:
\emph{Mordell-Weil lattices in characteristic 2. I. Construction and first properties},
Int.~Math.~Res.~Not.~1994, no.~8, pp.~343--361.



\bibitem{Gekeler} Gekeler, E.-U.: \emph{Local and global ramification
properties of elliptic curves in characteristics two and three}, in:
{\it Algorithmic Algebra and Number Theory} 
(Matzat, B. H., Greuel, G.-M., Hi\ss , G.\ eds.),
Berlin-Heidelberg-New York: Springer 1998, pp.~49--64


\bibitem{Hall}
Hall, M. Jr.:
\emph{The Diophantine equation $x\sp{3}-y\sp{2}=k$},  
Computers in number theory (Proc.~Sci.~Res.~Council Atlas Sympos.~No.~{\bf 2}, Oxford, 1969),  pp.~173--198.







%
%




\bibitem{Kodaira} Kodaira, K.: \emph{On compact analytic surfaces II, III},  Ann.~of Math.~(2) {\bf 77} (1963), 563--626; ibid.~{\bf 78} (1963), 1--40.






\bibitem{L} Livn\'e, R.: \emph{Motivic Orthogonal Two-dimensional Representations of 
Gal$(\bar{\Q}/\Q)$}, Israel J. Math. {\bf 92} (1995), pp.~149--156.













\bibitem{MP} Miranda, R., Persson, U.:
\emph{Configurations of $I_n$ Fibers on Elliptic K3 surfaces},
Math.~Z.~{\bf 201} (1989), pp.~339--361.



\bibitem{MT} Miyamoto, R., Top, J.: \emph{Reduction of Elliptic Curves in Equal Characteristic}, Canad.~Math.~Bull.~{\bf 48} (2005), pp.~428--444.



\bibitem{P-S} Pesenti, J., Szpiro, L.:  \emph{In\'egalit\'e du discriminant pour les pinceaux elliptiques \`a r\'eductions quelconques}, Compositio Math.~{\bf 120} (2000), no. 1, pp.~83--117.

\bibitem{PS}  Pjatecki\u\i -\v Sapiro, I.~I., \v Safarevi\v c, I.~R: \emph{A Torelli theorem for algebraic surfaces of type K3}, Math.~USSR Izv.~{\bf 5}, No.~3 (1972), pp.~547-588.





\bibitem{S-max} Sch\"utt, M.: \emph{The maximal singular fibres
of elliptic K3 surfaces}, Arch.~Math.~(Basel) {\bf 87}, Nr.~4 (2006), pp.~309--319.



\bibitem{S-CM} Sch\"utt, M.:
\emph{CM newforms with rational coefficients}, 
Ramanujan J.~{\bf 19} (2009), pp.~187--205. 




\bibitem{SS} Sch\"utt, M., Schweizer, A.: \emph{Davenport-Stothers inequalities and elliptic surfaces in positive characteristic}, Quarterly J.~Math.~{\bf 59} (2008), pp.~499--522.



\bibitem{ST} Sch\"utt, M., Top, J.: \emph{Arithmetic of the [19,1,1,1,1,1] fibration}, Comm.~Math.~Univ.~St.~Pauli {\bf 55}, 1 (2006), pp.~9--16.





\bibitem{Schw} Schweizer, A.: \emph{Extremal elliptic surfaces in 
characteristic $2$ and $3$,} Manuscripta Math. {\bf 102} (2000), pp.~505--521.





\bibitem{ShMW} Shioda, T.: \emph{On the Mordell-Weil lattices}, Comm.~Math.~Univ.~St.~Pauli, {\bf 39} (1990), pp.~211--240.



\bibitem{Sh-max} Shioda, T.: \emph{The elliptic K3 surfaces
with a maximal singular fibre,} C.~R.~Acad.~Sci. Paris,
Ser.~I, {\bf 337} (2003), pp.~461--466.





\bibitem{Sh-DS} Shioda, T.: \emph{Elliptic surfaces and Davenport-Stothers triples}, Comment.~Math.~Univ.~St.~Pauli {\bf 54} (2005), pp.~49--68.



\bibitem{SI} Shioda, T., Inose, H.:
\emph{On Singular $K3$ Surfaces},
in: Baily, W.~L.~Jr., Shioda, T.~(eds.),
\emph{Complex analysis and algebraic geometry},
Iwanami Shoten, Tokyo (1977), pp.~119--136.




\bibitem{Silverman} Silverman, J.~H.: \emph{Advanced Topics in the 
Arithmetic of Elliptic Curves}, Springer GTM, Berlin-Heidelberg-New York,
1994.



\bibitem{St} Stothers, W.~W.: \emph{Polynomial identities and
Hauptmoduln}, Quart.~J.~Math.~Oxford (2), {\bf 32} (1981), pp.~349--370.






\bibitem{Tate-C} Tate, J.: {\it Algebraic cycles and poles of zeta functions}, in: {\it Arithmetical Algebraic Geometry} (Proc. Conf. Purdue Univ., 1963), pp.~93--110, Harper \& Row (1965).








\bibitem{Tate} Tate, J.: {\it Algorithm for determining the type of a singular fibre in an elliptic pencil}, in: {\it Modular functions of one variable IV} (Antwerpen 1972), SLN {\bf 476} (1975), pp.~33--52.



\end{thebibliography}
\end{document}